\documentclass[11pt,a4paper]{article}
\usepackage{epsfig,psfrag,amsmath,amssymb,latexsym}
\usepackage{amscd,mathrsfs}
\usepackage{amsfonts}
\usepackage{graphicx}
\usepackage{epstopdf}
\usepackage{bbm}
\usepackage{sectsty}
\usepackage{amsthm}
\usepackage[numbers]{natbib}
\usepackage{dsfont}

\sectionfont{\large}

\def\P{\mathbb{P}}
\def\E{\mathbb{E}}

\def\EE{\mathcal{E}}
\def\N{\mathbb{N}}

\def\X{(X_n)_{n \in \N}}

\newcommand{\eop}{\mbox{ \vrule height7pt width7pt depth0pt}}

\newtheorem{definition}{Definition}[section]
\newtheorem{theorem}[definition]{Theorem}
\newtheorem{lemma}[definition]{Lemma}

\newtheorem{proposition}[definition]{Proposition}

\newtheorem{corollary}[definition]{Corollary}

\newtheorem*{pf}{Proof}
\newtheorem*{pf1}{Proof of Theorem \ref{orderHj}}
\newtheorem*{pf2}{Proof of Theorem \ref{spectralHs}}

\theoremstyle{remark}

\numberwithin{equation}{section}

\title{On hitting times for simple random walk on dense Erd\"os-R\'enyi random graphs}
\author{Matthias L\"owe and Felipe Torres \thanks{Corresponding author. E-mail address: ftorrestapia {\bf at} math.uni-muenster.de}
\\
{\small \it Institute for Mathematical Statistics, University of M\"unster} 
}
\begin{document}
\parindent 0em
\maketitle
\vspace{-12pt}
\abstract{\noindent
Let $(V,\EE)$ be a realization of the Erd\"os-R\'enyi random graph model $G(N,p)$ and $\X$ be a simple random walk on it. We study the size of $\sum_{i \in V} \pi_ih_{ij} $ where $\pi_i=d_i / 2|\EE|$ for $d_i$ the number of neighbors of node $i$ and $h_{ij}$ the hitting time for $\X$ between nodes $i$ and $j$. We always consider a regime of $p=p(N)$ such that realizations of $G(N,p)$ are asymptotically almost surely connected as $N\to\infty$. Our main result is that $\sum_{i \in V} \pi_ih_{ij} $ is almost surely of order $N(1+o(1))$ as $N\to \infty$. This coincides with previous non-rigorous results in the physics literature \cite{sood}. Our techniques are based on large deviations bounds on the number of neighbors of a typical node and the number of edges in $G(N,p)$ \cite{chungcomplex} together with bounds on the spectrum of the (random) adjacency matrix of $G(N,p)$ \cite{yau1}.
}
\section{Introduction} \label{introduction}
Random walks have been used since a couple of years to investigate properties of finite and infinite graphs e.g. \cite{doylesnell, lovaszgraphs, woess}, partially as part of a larger program for developing a theory of probability on finite and infinite graphs that accounts for their intrinsic geometrical structure, e.g. \cite{aldousfill, grimmett, levin, lyonsperes}, partially as an independent and exciting research area with its own rights. Meanwhile, the study of random graph models has received great attention not only within the probability community e.g. \cite{bollobas, chungcomplex, durret, jansonluczak, kolchinrandomgraphs, remco}, but also within the physics, biology, engineering, computer sciences and social sciences communities, among others. In this paper, we give a small contribution to that program by studying some properties of hitting times of a random walk on Erd\"os-R\'enyi random graphs.
\\
Let $\mathbb{G}:=(\mathcal{G},\mathcal{F}, \P_p)$ denote the probability space of the Erd\"os-R\'enyi random graph model $G(N,p)$ on $N$ vertices. More precisely, $\mathcal{G}$ is the set of all graphs on $N$ vertices, $\mathcal{F}$ is its powerset, and $\P_p$ is the probability measure for which every edge is created, independently one from each other, with probability $p \in [0,1]$. We will call a realization of such a graph $(V,\EE)$. We say that an event $A_N \subset \mathcal{G}$ happens asymptotically almost surely (abbreviated by a.a.s.) if $\P_p(A_N) \to 1$ as $N\to\infty$. Given $(V,\EE) \in \mathcal{G}$, let $\X$ be a simple random walk on $(V,\EE)$, i.e. $\X$ is the discrete time Markov Chain with state space $V$ and transition probabilities given by
\[ p^n_{ij}:=P(X_{n+1}=j \,|\,X_n=i)=\left\{ \begin{array}{ll}
1/d_i & \mbox{if $ij \in \EE$}\\
0 &\mbox{otherwise.}
\end{array}
\right.\]
Here $d_i$ is the degree of vertex $i\in V$. If $(V,\EE)$ is connected, it is well known that:
\begin{itemize}
\item $\X$ has a unique stationary distribution defined by the vector
$$\pi=(\pi_1,\dots,\pi_N)^T \quad \mbox{where  }\pi_i:=d_i / (2|\EE|).$$
\item For $i,j\in V$, $h_{ij}$ the expected number of steps $\X$ takes to visit $j$ when starting from $i$ is (a.s.) finite.
\end{itemize}
The quantity $h_{ij}$ is called {\it hitting (or access) time for $\X$ between $i$ and $j$}, see \eqref{tauj} below. Note that $h_{ij}$ is a random variable on $\mathbb{G}$ and is itself an expectation with respect to the law of $\X$ as well. In the present note, we would like to estimate the so called {\it random target time}, defined by
\begin{equation}\label{eh}
H_j:=\sum_{i \in V} \pi_ih_{ij}\,\,,
\end{equation}
for certain regime of $p=p(N)$ such that $(V,\EE)$ is a.a.s.\ connected. In \citep[chapter 10]{levin} connections between $H_j$ and some other relevant quantities (e.g.\ mixing time for random walks) are shown. Our main motivation is to give a rigorous proof of $H_j = N + o(N)$ for dense Erd\"os-R\'enyi graphs, as proposed in the physics work of \cite{sood}. In principle, $(V,\EE)$ is a.a.s.\ connected for $p=p_c:=\kappa \frac{\log N}{N}$ for $\kappa>1$ constant (cf. \cite{durret}, section 2.8 and references therein). In this note, it will be necessary to take $p(N)=\Omega ((\log N)^{2C\xi -1})$ (therefore much larger than $p_c$), because we use results on the spectrum of $A$ only available in this regime \cite{yau1} (cf. proof of Proposition \eqref{peg} in our Section 3). However, we conjecture that our results are true already for $p\ge p_c$.
\\
\\
The rest of this note is organized as follows. In Section \ref{defmainresults} we give the definitions and our main results on the order of magnitude of $H_j $, together with some applications estimating the order of magnitude of the so called {\it random starting time}, see \eqref{ehs}, and commenting on the order of magnitude of the commute time for $\X$, see \eqref{ct}. Section \ref{iresults} explains how to take advantage of an spectral decomposition for $h_{ij}$, for which every term will be bounded depending on the regime of $p=p(N)$ and finally the a.a.s.\ order of magnitude will be determined. The bounds on the degree of nodes and on the number of edges  are valid for every regime of $p=p(N)$, though they are sharper in the regime where $(V,\EE)$ becomes a.a.s.\ connected (and therefore also in the regime where our main theorem is stated for). Our technique could eventually be applied to the largest component in other regimes of $G(N,p)$, provided corresponding bounds on the spectrum of its (random) adjacency matrix.

\section{Definitions and main results}\label{defmainresults}
Let $(V,\EE)$ be an a.a.s.\ connected random graph in $\mathbb{G}$ and $\X$ be a simple random walk taking values on $V$ as defined above. Its law is denote by $\P$ and the corresponding expectation by $\E$, and as usual, for $i \in V$, let $\P_i(\cdot) =\P(\cdot | X_0=i)$ and $\E_i(\cdot)=\E(\cdot | X_0=i)$. As mentioned in the introduction, our aim is to estimate the order of magnitude of quantities involving hitting times of $\X$. Let $\tau_j$ be the first time $\X$ is at $j \in V$, i.e.
\begin{equation}
\tau_j:=\inf\{ m>0 : X_m=j\}.\label{tauj}
\end{equation}
The hitting time $h_{ij}$ is then given by $h_{ij}:=\E_i (\tau_j)$. We will represent $h_{ij}$ in terms of eigenvalues and eigenvectors of the adjacency matrix of $(V,\EE)$. Let $D$ be the diagonal matrix with entries $(D)_{ii}=1/d_i$ and $A$ be the adjacency matrix of $(V,\EE)$. Consider the symmetric (random) matrix 
$$
B:=D^{1/2} A D^{1/2}
$$ 
and its (random) eigenvalues $1=\lambda_1 \ge \lambda_2 \ge \cdots \ge \lambda_N$. The corresponding orthonormal (random) eigenvectors are denoted by $v_1,\dots,v_N$, and there components are $v_k=(v_{kj})_{j=1,\dots,N}$. Note that the positive (random) vector $w:=(\sqrt{d_1},\dots,\sqrt{d_N})$ satisfies $B \cdot w=1\cdot w$, therefore by the Perron-Frobenius Theorem $v_{1j}=\sqrt{d_j / 2|\EE|}$. Here $|\EE|$ is the number of edges in $(V,\EE)$. For later use, remark that:
\begin{itemize}
\item For every $k =2,\dots,N$
\begin{equation}\label{orthogonal}
0=v_k \cdot v_1 = \sum_{j=1}^N v_{kj}v_{1j}= \frac{1}{\sqrt{2|\EE|}}\sum_{j=1}^N v_{kj}\,\sqrt{d_j} \Rightarrow \sum_{j=1}^N v_{kj}\,\sqrt{d_j}=0
\end{equation}
\item The matrix $V:=(v_{kj})_{k,j=1,\dots,N}$ is unitary, hence its rows and its columns form an orthonormal set, therefore
\begin{equation}\label{upto1}
1= \sum_{j=1}^N v_{kj}^2 = \sum_{k=1}^N v_{kj}^2
\end{equation}
\end{itemize}
The following theorem gives a spectral decomposition of the hitting times:
\begin{theorem}{\rm \citep{lovaszgraphs}}\label{spectralH}
\[ h_{ij} = 2|\EE| \sum_{k=2}^N \frac{1}{1-\lambda_k} \left( \frac{v^2_{kj}}{d_j}-\frac{v_{ki}v_{kj}}{\sqrt{d_id_j}} \right),\]
a.a.s.\ for any $i,j \in V$, where $|\EE|$ is the number of edges in $G(N,p)$.
\end{theorem}

Our main results requires a quantity $\xi$ that stems from the estimation of the spectral gap of $A$: In what follows $a_0>0$ and $\tilde{a}_0 \ge 10$ will always be constants and $\xi=\xi_N$ will be a parameter such that
\begin{equation} \label{def xi}
1+a_0 \le \xi \le \tilde{a}_0\log\log N
\end{equation}
Our main result is:
\begin{theorem}\label{orderHj}
Let $\xi$ be as in \eqref{def xi}. Then, there exists a constant $C>0$ not depending on $N$ such that in the regime of $p=p(N)$ where $(\log N)^{2C\xi} / Np \to 0$ as $N \to \infty$, we have a.a.s.
\[H_j =N(1+o(1)).\]
\end{theorem}

Next, let us analyze quantities related to $H_j$. Given $i \in V$, let $H^i$ be the {\it random starting time} defined by
\begin{equation}\label{ehs}
H^i := \sum_{j\in V} \pi_j h_{ij}.
\end{equation}
Note that in general $h_{ij} \neq h_{ji}$ for $i,j \in V$, so therefore in general $H^i \neq H_i$. However, if the graph has a vertex-transitive automorphism group then $h_{ij}=h_{ji}$ for every $i,j \in V$ \citep[corollary 2.6]{lovaszgraphs}. From Theorem \ref{spectralH}  we can deduce that
\begin{eqnarray*}
H^i=\sum_{j=1}^N \pi_j h_{ij} &=& \sum_{j=1}^N \sum_{k=2}^N \frac{1}{1-\lambda_k} \left( v^2_{kj} - v_{kj}v_{ki}\sqrt{\frac{d_j}{d_i}}\,\right) \\
&=& \sum_{k=2}^N \frac{1}{1-\lambda_k} \left( \sum_{j=1}^N v^2_{kj} - v_{ki}\sqrt{\frac{1}{d_i}}\sum_{j=1}^N v_{kj} \sqrt{d_j}\right) \\
&=& \sum_{k=2}^N \frac{1}{1-\lambda_k}
\end{eqnarray*}
where the last equality follows from \eqref{orthogonal}. So we have for every $i \in V$,
\begin{equation}\label{sumH}
H^i = \sum_{k=2}^N \frac{1}{1-\lambda_k}.
\end{equation}
We therefore obtain the following order of magnitude for $H^i$:
\begin{theorem}\label{spectralHs}
Let $\xi$ be as in \eqref{def xi}. Then, there exists a constant $C>0$ not depending on $N$ such that in the regime of $p=p(N)$ where $(\log N)^{2C\xi} / Np \to 0$ as $N \to \infty$, we have a.a.s.
\[ H^i = N(1+o(1))\]
for each $i \in V$.
\end{theorem}
Another related quantity is the commute time. Let $\kappa_{ij}$ be the commute time between $i$ and $j$, i.e. the expected number of steps that $\X$, starting at $i \in V$, needs for coming back to $i$ visiting $j \in V$ before,
\begin{equation}\label{ct}
\kappa(i,j) := h_{ij}+h_{ji}.
\end{equation}
Due to the definition of the commute time in terms of hitting times and the spectral decomposition in Theorem \ref{spectralH}, we have \citep[corollary 3.2]{lovaszgraphs}
\[ \kappa(i,j) = 2|\EE| \sum_{k=2}^N \frac{1}{1-\lambda_k} \left( \frac{v_{kj}}{\sqrt{d_j}} - \frac{v_{ki}}{\sqrt{d_i}}\right)^2\]
and therefore we have \citep[corollary 3.3]{lovaszgraphs}
\begin{equation}\label{ctbounds}
|\EE| \left( \frac{1}{d_i} + \frac{1}{d_j} \right) \le \kappa(i,j) \le \frac{2|\EE|}{1-\lambda_2} \left( \frac{1}{d_i} + \frac{1}{d_j} \right).
\end{equation}
Remember that these results are also valid for random graphs though the graphs in \cite{lovaszgraphs} are deterministic. Therefore, combining Corollary \ref{c1pi}, Proposition \ref{peg} and inequality \eqref{ctbounds} we have:
\begin{corollary}\label{cto}
Let $\xi$ be as in \eqref{def xi}. Then, there exists a constant $C>0$ not depending on $N$ such that in the regime of $p=p(N)$ where $(\log N)^{2C\xi} / Np \to 0$ as $N \to \infty$, we have a.a.s.
\[ N(1+o(1)) \le \kappa(i,j) \le N(2+o(1))\]
for every $i,j \in V$.
\end{corollary}

\section{Intermediate results and proofs}\label{iresults}
\subsection{Spectral decomposition.}

For the readers' convenience, we give a brief summary over spectral decomposition of hitting times, to be found in \cite{lovaszgraphs}, which will be our starting point. Following the proof of Theorem 2.10 b) in \cite{lovaszgraphs}, we have
\begin{eqnarray*}
H_j = \sum_{i=1}^N \pi_i h_{ij} &=&\sum_{i=1}^N \sum_{k=2}^N \frac{1}{1-\lambda_k}\left( v^2_{kj}\frac{d_i}{d_j}-v_{ki}v_{kj}\sqrt{\frac{d_i}{d_j}}\,\right)\\
&=& \left(\,\frac{1}{d_j} \sum_{i=1}^N d_i\,\right)\sum_{k=2}^N \frac{1}{1-\lambda_k} v^2_{kj} - \sum_{k=2}^N v_{kj} \frac{1}{\sqrt{d_j}} \sum_{i=1}^Nv_{ki} \sqrt{d_i}\\
&=& \left(\,\frac{1}{d_j} \sum_{i=1}^N d_i\,\right)\sum_{k=2}^N \frac{1}{1-\lambda_k} v^2_{kj}
\end{eqnarray*}
where the last inequality follows because of \eqref{orthogonal}. At last,
\begin{equation}\label{help0}
H_j = \frac{2|\EE|}{d_j} \sum_{k=2}^N \frac{1}{1-\lambda_k}v^2_{kj}.
\end{equation}
Now, note that by using \eqref{upto1} we can re-write
\begin{equation}\label{help1}
\sum_{k=2}^N v^2_{kj} = \sum_{k=1}^N v^2_{kj} - \pi_j = 1-\pi_j
\end{equation}
\begin{equation}\label{help2}
\sum_{k=2}^N (1-\lambda_k) v^2_{kj} = \sum_{k=1}^N (1-\lambda_k) v^2_{kj} = 1-\sum_{k=1}^N \lambda_k v^2_{kj} = 1-(B)_{jj}=1.
\end{equation}
Therefore, from the inequality between the arithmetic and harmonic means (considering $v^2_{kj}$ as weights) which can be written as
\[ \frac{\sum_{k=2}^N \frac{1}{1-\lambda_k} v^2_{kj}}{\sum_{k=2}^N v^2_{kj}} \ge \frac{\sum_{k=2}^N v^2_{kj}}{\sum_{k=2}^N (1-\lambda_k)v^2_{kj}}\]
and by using \eqref{help1}, \eqref{help2} and \eqref{help0} we have that
\[H_j = \frac{1}{\pi_j}\sum_{k=2}^N \frac{1}{1-\lambda_k} v^2_{kj} \ge \frac{1}{\pi_j}(1-\pi_j)^2 = \frac{1}{\pi_j} -2 + \pi_j \ge \frac{1}{\pi_j} -2=\frac{2|\EE|}{d_j}-2,\]
which is a lower bound for $H_j$ only in terms of $|\EE|$ and $d_j$, namely
\begin{equation}\label{lowerHj}
H_j \ge \frac{2|\EE|}{d_j} -2.
\end{equation}
An upper bound for $H_j$ can be obtained by starting from \eqref{help0} as follows:
\begin{eqnarray*}
H_j &=& \frac{2|\EE|}{d_j} \sum_{k=2}^N \frac{1}{1-\lambda_k}v^2_{kj} \\
&\le& \frac{2|\EE|}{d_j} \cdot\frac{1}{1-|\lambda_2|} \sum_{k=2}^N v^2_{kj} \\
&=& \frac{2|\EE|}{d_j} \cdot\frac{1}{1-|\lambda_2|} \left( \sum_{k=1}^N v^2_{kj} - \pi_j\right) = \frac{2|\EE|}{d_j} \cdot\frac{1}{1-|\lambda_2|} \left(1-\frac{d_j}{2|\EE|}\right)
\end{eqnarray*}
which can be written as
\begin{equation}\label{upperHj}
H_j \le \left(\frac{2|\EE|}{d_j}-1 \right) \left(\frac{1}{1-|\lambda_2|}\right).
\end{equation}
Hence, \eqref{lowerHj} and \eqref{upperHj} show that:
\begin{proposition}\label{boundsHj}
For every $j \in V$, we have
\[ \frac{2|\EE|}{d_j} -2\le H_j \le  \left(\frac{2|\EE|}{d_j}-1 \right) \left(\frac{1}{1-|\lambda_2|}\right).\]
\end{proposition}
The next step is now obvious: We will estimate the almost sure order of magnitude of $2|\EE| / d_j$ and $(1-|\lambda_2|)^{-1}$ (the inverse of the spectral gap) as $N\to\infty$.
\subsection{Order of the stationary distribution.}

The following two results are valid for all regimes of $p=p(N)$, and even hold for a more general model of weighted random graphs. The first one, is about large deviations for the random variable $d_j$:
\begin{lemma}{\rm \citep{chungcomplex}}\label{lemma1}
For $(V,\EE)$, with probability $1-e^{-c^2/2}$, the degree $d_j$ satisfies
\begin{equation}\label{lowerd}
d_j > Np - c\sqrt{Np},
\end{equation}
and with probability $1-e^{-\frac{c^2}{2(1+c/(3\sqrt{Np}))}}$ the degree $d_j$ satisfies
\begin{equation}\label{upperd}
d_j < Np + c\sqrt{Np}.
\end{equation}
\end{lemma}
Lemma \ref{lemma1} provides useful inequalities in the regime $N p > c' \log N$ with $c'>1$ constant. The second result is about large deviations for the number of edges $|\EE|$:
\begin{lemma}{\rm \citep{chungcomplex}}\label{lemma2}
With probability $1-2e^{-c^2/6}$, the number $|\EE|$ of edges in $(V,\EE)$ satisfies
\begin{equation}\label{boundE}
\Big| 2|\EE| - N^2p \Big| < c\sqrt{N^2 p}
\end{equation}
for every $0<c\le\sqrt{N^2 p}$.
\end{lemma}
As a consequence, we can deduce:
\begin{corollary}\label{c1pi}
In the regime of $p=p(N)$ such that $\frac{\log N}{Np} \to 0$ as $N \to \infty$, we have a.a.s.
\begin{equation*}\label{stat}
\frac{2|\EE|}{d_j} = N (1+o(1))
\end{equation*}
for every $j\in V$.
\end{corollary}
\begin{pf}
Take $c=\sqrt{\log N}$ in \eqref{lowerd} and in \eqref{boundE}. Then, we have that a.a.s.
\[\frac{2|\EE|}{d_j} \le \frac{N^2 p + \sqrt{N^2p\log N}}{Np - \sqrt{Np\log N}} = N \left( \frac{1+\sqrt{\frac{\log N}{Np}}\frac{1}{\sqrt N}}{1-\sqrt{\frac{\log N}{Np}}}\right)=N(1+o(1)).\]
For the other inequality, again taking $c=\sqrt{\log N}$, now in \eqref{upperd} and in \eqref{boundE}, we have a.a.s.\
\[\frac{2|\EE|}{d_j} \ge \frac{N^2 p - \sqrt{N^2p\log N}}{Np + \sqrt{Np\log N}} = N \left( \frac{1-\sqrt{\frac{\log N}{Np}}\frac{1}{\sqrt N}}{1+\sqrt{\frac{\log N}{Np}}}\right)=N(1+o(1)).\quad\eop\]
\end{pf}
\subsection{Order of the spectral gap.}

The idea here is to use the a recent analysis on eigenvalues statistics of the adjacency matrix of dense Erd\"os-R\'enyi graphs. This was obtained in \cite{yau1} to estimate the order of $(1-|\lambda_2|)^{-1}$. Let us start by deriving an a.a.s.\ relation between the eigenvalues of the matrix $A=(a_{ij})_{i,j=1,\dots,N}$ and the eigenvalues of the matrix $B=D^{1/2}AD^{1/2}=(b_{ij})_{i,j=1,\dots,N}$. Note that
\[ b_{ij}=\frac{a_{ij}}{\sqrt{d_i d_j}}=\frac{1}{Np}a_{ij} + \left(\frac{Np-\sqrt{d_id_j}}{Np\sqrt{d_id_j}}\right)a_{ij}.\]
Therefore, $B = \frac{1}{Np}A + R$ where the matrix $R=(r_{ij})_{i,j=1,\dots,N}$ is defined by
\[ r_{ij}=\left(\frac{Np-\sqrt{d_id_j}}{Np\sqrt{d_id_j}}\right)a_{ij}.\]
Recall that if $M=(m_{ij})$ is an $N \times N$-matrix with real entries, the spectral radius of $M$ given by
$$\rho(M):=\max\{ | \lambda | : \mbox{$\lambda$ eigenvalue of $M$} \}$$
and it holds $\rho(M) \le || M ||_\infty:= \max_{1 \le i \le N} \sum_{j=1}^N m_{ij} $. Let $\nu_1 \ge \nu_2 \ge \cdots \ge \nu_N$ be the eigenvalues of $A$, and $w_1,w_2,\dots,w_N$ be their corresponding normalized eigenvectors. We are interested in the case where $(V,\EE)$ is (a.a.s.) not the complete graph, such that $A$ will not have $(1/N,\dots,1/N)$ as normalized eigenvector. Now, it is easy to see that
\[ w_2^T Bw_2 = \frac{\nu_2}{Np} + w_2^TRw_2\]
therefore
\begin{eqnarray*}
|\lambda_2| =\max_{w_2 : {\footnotesize \begin{array}{c}w_2^Tw_2=1 \\ w_2 \bot v_1 \end{array}}} |w_2^T B w_2| &\le& \frac{|\nu_2|}{Np} +  \max_{w_2 : w_2^Tw_2=1 } |w_2^T R w_2| \\
&=& \frac{|\nu_2|}{Np} + \rho(R) \le \frac{|\nu_2|}{Np} + ||R||_\infty
\end{eqnarray*}
which means that
\begin{equation}\label{upperB}
|\lambda_2| \le \frac{|\nu_2|}{Np} +\max_{i=1,\dots,N}\sum_{j=1}^N \frac{|Np-\sqrt{d_id_j}|}{Np\sqrt{d_id_j}}a_{ij}.
\end{equation}
\begin{lemma}\label{Bspectral}
In the regime of $p=p(N)$ such that $\frac{\log N}{Np} \to 0$ as $N \to \infty$, we have that a.a.s.\ $\lambda_2$ and $\nu_2$ satisfy
\[ |\lambda_2| \le \frac{|\nu_2|}{Np} + o(1).\]
\end{lemma}
\begin{pf}
Because of inequality \eqref{upperB}, it is enough to show that a.a.s.\ we have
\[ \max_{i=1,\dots,N}\sum_{j=1}^N \frac{|Np-\sqrt{d_id_j}|}{Np\sqrt{d_id_j}}a_{ij} =o(1)\]
for any $j\in V$. Indeed, setting $c:=\sqrt{\log N}$, from Lemma \ref{lemma1} we get that
\[| Np - \sqrt{d_id_j}| \le \sqrt{Np\log N}\]
with probability $(1-\frac{1}{\sqrt{N}})(1-2\exp(-\frac{1}{4}(\frac{1}{\log N}+\sqrt{\frac{1}{9Np \log N}})^{-1}))$ which tends to 1 as $N\to\infty$. Moreover, we obtain
\[ d_i d_j \ge (Np - \sqrt{Np\log N})^2\]
with probability $(1-\frac{1}{\sqrt{N}})^2$ which tends to 1 as $N\to\infty$. By combining these two inequalities, given any $\varepsilon>0$, there exists $N_0=N_0(\varepsilon)>0$ such that
\begin{eqnarray*}
\sum_{j=1}^N \frac{|Np-\sqrt{d_id_j}|}{Np\sqrt{d_id_j}}a_{ij} &\le& \frac{\sqrt{Np\log N}}{Np(Np-\sqrt{Np\log N})}\sum_{j=1}^Na_{ij} \\
&=& \frac{\sqrt{Np\log N}}{Np(Np-\sqrt{Np\log N})}d_i \\
&\le& \frac{\sqrt{Np\log N}(Np +\sqrt{Np\log N})}{Np(Np-\sqrt{Np\log N})}\\
&=&\sqrt{\frac{\log N}{Np}} \cdot \left( \frac{1+\sqrt{\frac{\log N}{Np}}}{1-\sqrt{\frac{\log N}{Np}}}\right) < \varepsilon
\end{eqnarray*}
a.s.\ for every $N > N_0$.\eop
\end{pf}
The next result establishes a bound for $(1-|\lambda_2|)^{-1}$:
\begin{proposition}\label{peg}
In the regime of $p=p(N)$ such that $(\log N)^{2C\xi} / Np \to 0$ as $N \to \infty$, we have that a.s.s.
\begin{equation*}\label{spectralgap}
\frac{1}{1-|\lambda_2|} \le 1+ o(1).
\end{equation*}
\end{proposition}
\begin{pf}
Let $\tilde{A}:=(Np(1-p))^{-1/2} A$. Let $\mu_1 \le \cdots \le \mu_N$ be the eigenvalues of $\tilde{A}$. One checks that the matrix $\tilde{A}$ satisfies Definition 2.2 in \cite{yau1} with $f=\sqrt{Np/(1-p)}$. On the other hand,  it is clear that
\[ \tilde{A} w_2 = (Np(1-p))^{-1/2} A w_2= (Np(1-p))^{-1/2} \nu_2 w_2\]
therefore $\mu_{N-1} = (Np(1-p))^{-1/2} \nu_2$. Moreover, inequality (3.19) in \cite{yau2} tells us that there exist constants $\theta,C>0$ (both not depending on $N$) such that
\begin{equation}\label{mumax}
|\mu_{N-1}| \le 2 + (\log N)^{C\xi} \left(\frac{1}{Np} + \frac{1}{N^{2/3}}\right)
\end{equation}
with probability at least $1-e^{-\theta (\log N)^\xi}$, which tends to one as $N \to \infty$. Combining this with \eqref{mumax} and Lemma \ref{Bspectral}, we have that
\begin{eqnarray*}
|\lambda_2| &\le& \sqrt{\frac{1-p}{Np}} \left[2 + (\log N)^{C\xi} \left(\frac{1}{Np} + \frac{1}{N^{2/3}}\right)\right] + o(1) \\
&\le& \sqrt{\frac{4}{Np}} + \sqrt{\frac{(\log N)^{2C\xi}}{Np}}\left(\frac{1}{Np} + \frac{1}{N^{2/3}}\right)+ o(1) \\
&=& o(1)
\end{eqnarray*}
a.s.s.\ where the last $o(1)$ is due to the hypothesis $(\log N)^{2C\xi} / Np \to 0$ as $N \to \infty$. So, a.a.s.\ we have that
\[1-|\lambda_2| \ge 1-o(1) \Rightarrow \frac{1}{1-|\lambda_2|} \le \frac{1}{1-o(1)} = 1+\left( \frac{1}{1-o(1)} -1\right)=1+o(1).\quad\eop\]
\end{pf}
\subsection{Proof of the main results.}

\begin{pf1}
Take $j \in V$. The regime in Proposition \ref{peg} also satisfies the condition $\frac{\log N}{Np} \to 0$ as $N \to \infty$, so a.a.s. we have that
\[H_j \le  \left(\frac{2|\EE|}{d_j}-1 \right) \left(\frac{1}{1-|\lambda_2|}\right) \le [N(1+o(1))-1](1+o(1)) = N(1+o(1))\]
which follows from the upper bound in Proposition \ref{boundsHj} and Proposition \ref{peg}, and
\[ N(1+o(1)) = N(1+o(1))-2 =\frac{2|\EE|}{d_j} -2 \le H_j\]
which follows from the lower bound in Proposition \ref{boundsHj} and Corollary \ref{c1pi}.\eop
\end{pf1}
\begin{pf2} In this regime, we use Proposition \ref{peg} and \eqref{sumH} to deduce a.a.s.\ the upper bound
\[ \sum_{k=2}^N \frac{1}{1-\lambda_k} \le \frac{1}{1-|\lambda_2|}(N-1) =(1+o(1))(N-1) = N(1+o(1)).\]
For the lower bound, note that $\frac 1{1-\lambda_k}1_{\{\lambda_k\ge 0\}}$ is positive, and under $1_{\{\lambda_k\ge 0\}}$ we are in the case $0 \le \lambda_k < 1$ for $k=2,\dots,N$. Therefore $1 \ge 1-\lambda_k >0$ and so it follows $\frac 1{1-\lambda_k}1_{\{\lambda_k\ge 0\}} \ge 1$, which implies
\[ \sum_{k=2}^N \frac{1}{1-\lambda_k} \ge \sum_{k=2}^N \frac{1}{1-\lambda_k} {\bf 1}_{\{\lambda_k \ge 0\}} \ge (N-1) = N(1+o(1)).\qquad \eop\]
\end{pf2}


\section*{Acknowledgements}
Felipe Torres would like to thank the support of the DFG (Deutsche Forschungsgemeinschaft) through the Research Cooperative Center (SFB) 878 "Groups, Geometry \& Actions" at University of M\"unster, Germany.


{\footnotesize \bibliographystyle{alpha}			
\bibliography{AT}}

\begin{thebibliography}{EKYY12}

\bibitem[AF00]{aldousfill}
David Aldous and Jim Fill.
\newblock Reversible markov chains and random walks on graphs. book in
  preparation.
\newblock {\em URL for draft at http://www. stat. Berkeley. edu/users/aldous},
  2000.

\bibitem[Bol01]{bollobas}
B{\'e}la Bollob{\'a}s.
\newblock {\em Random graphs}, volume~73.
\newblock Cambridge university press, 2001.

\bibitem[CL06]{chungcomplex}
Fan Chung and Linyuan Lu.
\newblock {\em Complex graphs and networks (cbms regional conference series in
  mathematics)}.
\newblock American Mathematical Society, Boston, MA, 2006.

\bibitem[DS84]{doylesnell}
Peter~G Doyle and J~Laurie Snell.
\newblock {\em Random walks and electric networks}, volume~22.
\newblock Mathematical association of America, 1984.

\bibitem[Dur06]{durret}
Rick Durrett.
\newblock {\em Random graph dynamics}, volume~20.
\newblock Cambridge University Press, 2006.

\bibitem[EKYY11]{yau1}
L{\'a}szl{\'o} Erd{\H{o}}s, Antti Knowles, Horng-Tzer Yau, and Jun Yin.
\newblock Spectral statistics of erd{\H{o}}s-r{\'e}nyi graphs i: Local
  semicircle law.
\newblock {\em arXiv preprint arXiv:1103.1919}, 2011.

\bibitem[EKYY12]{yau2}
L{\'a}szl{\'o} Erd{\H{o}}s, Antti Knowles, Horng-Tzer Yau, and Jun Yin.
\newblock Spectral statistics of erd{\H{o}}s-r{\'e}nyi graphs ii: Eigenvalue
  spacing and the extreme eigenvalues.
\newblock {\em Communications in Mathematical Physics}, pages 1--54, 2012.

\bibitem[ER59]{erdosrenyiFIRST}
Paul Erd{\H{o}}s and Alfr{\'e}d R{\'e}nyi.
\newblock On random graphs.
\newblock {\em Publicationes Mathematicae Debrecen}, 6:290--297, 1959.

\bibitem[Gri10]{grimmett}
Geoffrey Grimmett.
\newblock {\em Probability on graphs: random processes on graphs and lattices},
  volume~1.
\newblock Cambridge University Press, 2010.

\bibitem[J{\L}R00]{jansonluczak}
Svante Janson, Tomasz {\L}uczak, and Andrzej Ruci{\'n}ski.
\newblock Random graphs. 2000.
\newblock {\em Wiley--Intersci. Ser. Discrete Math. Optim}, 2000.

\bibitem[Kol99]{kolchinrandomgraphs}
Valentin~Fedorovich Kolchin.
\newblock {\em Random Graphs, Encyclopedia of Mathematics and its Applications,
  53}.
\newblock Cambridge University Press, Cambridge, 1999.

\bibitem[Lov93]{lovaszgraphs}
Laszlo Lov{\'a}sz.
\newblock Random walks on graphs: A survey.
\newblock {\em Combinatorics, Paul Erd\"os is eighty}, 2(1):1--46, 1993.

\bibitem[LP09]{lyonsperes}
Russell Lyons and Yuval Peres.
\newblock Probability on trees and networks.
\newblock {\em In preparation. Current version available at http://mypage. iu.
  edu/\~{} rdlyons/prbtree/book. pdf}, 2009.

\bibitem[LPW09a]{levinperes}
David~Asher Levin, Yuval Peres, and Elizabeth~Lee Wilmer.
\newblock {\em Markov chains and mixing times}.
\newblock AMS Bookstore, 2009.

\bibitem[LPW09b]{levin}
David~Asher Levin, Yuval Peres, and Elizabeth~Lee Wilmer.
\newblock {\em Markov chains and mixing times}.
\newblock Amer Mathematical Society, 2009.

\bibitem[LT13]{loewetorresscGWt}
Matthias L\"owe and Felipe Torres.
\newblock A note on hitting times for simple random walk on rooted, subcritical
  critical galton-watson trees.
\newblock {\em In preparation}, October, 2013.

\bibitem[LT14]{loewetorresscER}
Matthias L\"owe and Felipe Torres.
\newblock On hitting times for simple random walk on critical and supercritical
  erd\"os-r\'enyi random graphs.
\newblock {\em In preparation}, January, 2014.

\bibitem[SRBA04]{sood}
Vishal Sood, Sidney Redner, and Dani Ben-Avraham.
\newblock First-passage properties of the erd{\H{o}}s--renyi random graph.
\newblock {\em Journal of Physics A: Mathematical and General}, 38(1):109,
  2004.

\bibitem[VDH09]{remco}
Remco Van Der~Hofstad.
\newblock Random graphs and complex networks.
\newblock {\em Available at http://www. win. tue. nl/rhofstad/NotesRGCN. pdf},
  2009.

\bibitem[Woe00]{woess}
Wolfgang Woess.
\newblock {\em Random walks on infinite graphs and groups}, volume 138.
\newblock Cambridge university press, 2000.

\end{thebibliography}

\end{document}